\documentclass[preprint,12pt]{elsarticle}

\usepackage{amssymb}
\newcommand{\la}{\lambda}
\newcommand{\vp}{{\mathbf p}}
\newcommand{\vz}{{\mathbf z}}
\newcommand{\vf}{{\mathbf f}}
\newcommand{\vx}{{\mathbf x}}

\newcounter{theorem}

\newcounter{theoremcounter}

\newcounter{remarkcounter}

\newcounter{examplecounter}

\newtheorem{theorem}[theoremcounter]{Theorem}

\newtheorem{remark}[remarkcounter]{Remark}

\newtheorem{example}[examplecounter]{Example}

\begin{document}

\begin{frontmatter}

%% Title, authors and addresses

%% use the tnoteref command within \title for footnotes;
%% use the tnotetext command for the associated footnote;
%% use the fnref command within \author or \address for footnotes;
%% use the fntext command for the associated footnote;
%% use the corref command within \author for corresponding author footnotes;
%% use the cortext command for the associated footnote;
%% use the ead command for the email address,
%% and the form \ead[url] for the home page:
%%
%% \title{Title\tnoteref{label1}}
%% \tnotetext[label1]{}
%% \author{Name\corref{cor1}\fnref{label2}}
%% \ead{email address}
%% \ead[url]{home page}
%% \fntext[label2]{}
%% \cortext[cor1]{}
%% \address{Address\fnref{label3}}
%% \fntext[label3]{}

\title{Two-sided bounds on the rate of convergence  for continuous-time finite inhomogeneous Markov chains}

%% use optional labels to link authors explicitly to addresses:
%% \author[label1,label2]{<author name>}
%% \address[label1]{<address>}
%% \address[label2]{<address>}
\author[zai]{A. I. Zeifman}
\author[vyk]{V. Yu. Korolev}
\address[zai]{Vologda State University;
Institute of Informatics Problems, Russian Academy of Sciences; Institute of Socio-Economic
Development of Territories, Russian Academy of Sciences; e-mail a$\_$zeifman@mail.ru}
\address[vyk]{Faculty of Computational Mathematics and
Cybernetics, Lomonosov Moscow State University; Institute of
Informatics Problems, Russian Academy of Sciences}

\begin{abstract}
%% Text of abstract
We suggest an approach to obtaining general two-sided bounds on the
rate of convergence  in terms of special ``weighted'' norms related
to total variation. Some important classes of continuous-time Markov
chains are considered: birth-death-catastrophes processes, queueing
models with batch arrivals and group services, chains with
absorption in zero.
\end{abstract}

\begin{keyword}
%% keywords here, in the form: keyword \sep keyword
continuous-time Markov chains \sep inhomogeneous Markov chains \sep
ergodicity bounds \sep special norms
%% MSC codes here, in the form: \MSC code \sep code
%% or \MSC[2008] code \sep code (2000 is the default)

\MSC[2010] 60J271 \sep 60J22 \sep 60J28
\end{keyword}

\end{frontmatter}

%%
%% Start line numbering here if you want
%%
% \linenumbers

%% main text
\section{Introduction}

The problem of finding sharp bounds on the rate of convergence to
the limiting characteristics for Markov chains is very  important
for at least two following reasons:

(i) Beginning from what point can be the exact characteristics of the
process be replaced by the limiting characteristics that are much easier
to calculate? To find the bounds on this time moment we need sharp upper (or both upper and
lower exact) bounds on the rate of convergence.

(ii) Perturbation bounds play a significant role in applications.
It is well known that even for homogeneous chains the best bounds
require the corresponding best bounds on the rate of convergence
(Kartashov, 1985, 1996;  Liu, 2012; Mitrophanov, 2003, 2004). Bounds
for general inhomogeneous Markov chains are also based on the estimates
of the rate of convergence in the special weighted norms (Zeifman and Korolev, 2014a).

In this note we deal with the weak ergodicity of finite inhomogeneous
continuous-time Markov chains and formulate an algorithm for
obtaining sharp upper and lower bounds on the rate of convergence in
the special weighted norms related to total variation. Our general
approach is closely connected with the notion of logarithmic norm and the
corresponding bounds for the Cauchy matrix (Zeifman, 1985, 1995;
Granovsky and Zeifman, 1997, 2004;  Zeifman et al, 2006, 2013; 2014b), and
allows to obtain general upper and lower bounds on the rate of
convergence for a countable state space and general initial conditions. On
the other hand, if the state space of the chain is finite, then the desired
bounds can be obtained in a simple way without this special technique.

Section 2 contains preliminary material. The general result is
presented is Section 3. In Sections 4-6 three important classes of
continuous-time Markov chains are considered for which it is
possible to obtain exact convergence rate estimates.

\section{Preliminaries}

Let $X=X(t)$, $t\geq 0$ be an inhomogeneous finite
continuous-time Markov chain. Let $p_{ij}(s,t)=\Pr\left\{
X(t)=j\left| X(s)=i\right. \right\}$, $0 \le i,j \le S$, $0\leq
s\leq t$, be the transition probabilities for $X=X(t)$,
$p_i(t)=\Pr\left\{ X(t) =i \right\}$  be its state probabilities,
and ${\bf p}(t) = \left(p_0(t), p_1(t), \dots, p_S(t)\right)^T$ be
the corresponding probability distribution. Throughout the paper we
assume that
\begin{equation}
\Pr\left\{X\left( t+h\right) =j|X\left( t\right) =i\right\} =\!
\begin{cases}{q_{ij}\left( t\right) h+\alpha_{ij}\left(t, h\right), & $\mbox{ if
}j\neq i$,\vspace{2mm}\cr
\\ 1\!-\!\sum\limits_{k\neq i}q_{ik}\left( t\right)
h+\alpha_{i}\left( t,h\right), & $\mbox{ if } j=i$,} \end{cases}
\label{1001}
\end{equation}
\noindent where $q_{ij}(t)$ denotes the intensity of transition from
the state $i$ to the state $j$, $\alpha_{ij}(t,h) = o(h)$ for any $t$ as $h
\to 0$. We assume that (in the inhomogeneous case) all intensity functions
are locally integrable on $[0,\infty)$.

Let $a_{ij}(t) =  q_{ji}(t)$ for $j\neq i$ and let $a_{ii}(t) =
-\sum_{j\neq i} a_{ji}(t) = -\sum_{j\neq i} q_{ij}(t)$.  The
probabilistic dynamics of the process is represented by the forward
Kolmogorov system
\begin{equation} \label{ur01}
\frac{d\vp}{dt}=A(t)\vp(t),
\end{equation}
\noindent where $A(t)$ is the transposed intensity matrix of the
process.

Throughout the paper by $\|\cdot\|$  we denote  the $l_1$-norm, i.e.
$\|{\vx}\|=\sum|x_i|$, and $\|B\| = \sup_j \sum_i |b_{ij}|$ for $B =
(b_{ij})_{i,j=0}^{n}$. Let $\Omega$ be the set all stochastic
vectors, i.e., vectors with nonnegative coordinates and unit
$l_1$-norm. It is well known  that the Cauchy problem for
differential equation (\ref{ur01}) has a unique solution for an
arbitrary initial condition, and  $\vp(s) \in \Omega$ implies
$\vp(t) \in \Omega$ for $t \ge s \ge 0$. Hence, we can put $p_0(t) =
1 - \sum_{i \ge 1} p_i(t)$, and obtain from (\ref{ur01}) the
equation (see detailed discussion in (Granovsky and Zeifman, 2004;
Zeifman et al, 2006))
\begin{equation}
\frac{d\vz}{dt}= B(t)\vz(t)+\vf(t), \label{2.06}
\end{equation}
\noindent where $\vf(t)=\left(a_{10}, a_{20},\cdots, a_{S0}
\right)^T$, $\vz(t)=\left(p_{1}, p_{2},\cdots, p_{S} \right)^T$,
\begin{equation}
{%\scriptsize
B = \left(
\begin{array}{ccccccccc}
a_{11}- a_{10}   & a_{12} - a_{10}   &  \cdots & a_{1S} - a_{10}  \\
a_{21} - a_{20} & a_{22} - a_{20}   &    \cdots & a_{2S} - a_{20}  \\
a_{31} - a_{30}    & a_{32} - a_{30}  &    \cdots & a_{3S} - a_{30}  \\
\cdots \\
a_{S1} - a_{S0}  & a_{S2} - a_{S0} & \cdots     &  a_{SS} -a_{S0}
\end{array}
\right).}\label{2.07}
\end{equation}

Recall that an inhomogeneous Markov chain $X(t)$ is called {\bf
weakly ergodic}, if $\|{\bf p^{*}}(t) - {\bf p^{**}}(t)\| \to 0$ as
$t \to \infty$  for any initial conditions ${\bf p^{*}}(0)$, ${\bf
p^{**}}(0)$, where ${\bf p^{*}}$, ${\bf p^{**}}$ are two
distributions of state probabilities for $X(t)$. A matrix $H$ is called
{\bf essentially nonnegative} if all its off-diagonal elements are nonnegative.

\section{General bounds}

Let  $D$ be a nonnegative matrix  such that

\noindent $(i)$ $\|{\bf z}\|_{1D} = \|D{\bf z}\| \ge d  \|{\bf z}\|$
for some positive $d$ and

\noindent $(ii)$ $DB(t)D^{-1}$ is essentially nonnegative for any $t \ge 0$.

It is worth noting that
$${\bf p}(t)= \left(
\begin{array}{cc} 1-
\sum_{i=1}^{S}z_i(t) \\ {\bf z}(t)\end{array} \right),$$ and hence
${\bf z}(t)$ completely determines the state probabilities
for $X(t)$.

Let $DB(t)D^{-1} = H(t) = \left(h_{ij}(t)\right)_{i,j=1}^S$.
Consider the quantities
\begin{equation}
h^*(t)=\max_j\sum_i r_{ij}, \quad  h_*(t)=\min_j\sum_i r_{ij}
\label{301}
\end{equation}

 \begin{theorem} Let $X(t)$ be a given finite inhomogeneous Markov chain,
and let there exist a matrix $D$ satisfying $(i)$--$(ii)$.
Then the following bounds hold:
\begin{equation}
\|{\bf z^{*}}(t) - {\bf z^{**}}(t)\|_{1D} \le e^{\int_0^t
h^*(\tau)\,d\tau}\|{\bf z^{*}}(0) - {\bf z^{**}}(0)\|_{1D},
\label{302}
\end{equation}
\noindent for any corresponding initial conditions ${\bf p^{*}}(0)$,
${\bf p^{**}}(0)$, and
\begin{equation}
\|{\bf z^{*}}(t) - {\bf z^{**}}(t)\|_{1D} \ge e^{\int_0^t
h_*(\tau)\,d\tau}\|{\bf z^{*}}(0) - {\bf z^{**}}(0)\|_{1D},
\label{303}
\end{equation}
\noindent if the initial conditions are such that $D\left({\bf z^{*}}(0)
- {\bf z^{**}}(0)\right) \ge {\bf 0}$.
\end{theorem}

{\bf Proof.} Put ${\bf v}(t)={\bf z^{*}}(t) - {\bf z^{**}}(t)$, and
let ${\bf x}(t) = D {\bf v}(t).$ Then
\begin{equation}
\frac{d{\bf x}}{dt}=H(t){\bf x}. \label{304}
\end{equation}
First, let ${\bf x}(0) \ge {\bf 0}$. Then ${\bf x}(t) \ge {\bf 0}$
for any $t \ge 0$, hence $\|{\bf x}(t)\| = \sum_i x_i(t)$. From (\ref{304}) we
obtain that
$$\frac{d\|{\bf x}\|}{dt} = \frac{d{\sum_i x_i}}{dt} = \sum_{i}\Big(\sum_{j} h_{ij}x_{j}\Big) =
\sum_{j}\Big(\sum_{i} h_{ij}\Big)x_{j} \le h^*(t) \sum_{j} x_j = h^*(t) \|{\bf x}\|.$$
Then $\|{\bf x}(t)\| \le e^{\int_0^t h^*(\tau)\,d\tau}\|{\bf
x}(0)\|$, if ${\bf x}(0) \ge {\bf 0}$. Let now ${\bf x}(0)$ be
arbitrary. Put $x^+_i(0) = \max \left(x_i(0), 0\right)$, ${\bf
x^+}(0) = \left(x^+_1(0),\dots,x^+_n(0)\right)^T$, and ${\bf x^-}(0)
= {\bf x}(0) - {\bf x^+}(0)$. Then ${\bf x^+}(0)  \ge {\bf 0},$
${\bf x^-}(0) \ge 0$, ${\bf x}(0) = {\bf x^+}(0) - {\bf x^-}(0)$,
and, moreover, $\|{\bf x}(0)\| = \|{\bf x^+}(0)\|+ \|{\bf x^-}(0)\|$.
Therefore,
$$\|{\bf x}(t)\| = \|{\bf x^+}(t) - {\bf x^-}(t)\| \le \|{\bf x^+}(t)\|+
\|{\bf x^-}(t)\| \le $$
$$e^{\int_0^t h^*(\tau)\,d\tau}\left(\|{\bf x^+}(0)\|+ \|{\bf x^-}(0)\|\right) =e^{\int_0^t h^*(\tau)\,d\tau}\|{\bf x}(0)\|.$$
On the other hand, if ${\bf x}(0) \ge {\bf 0}$, then $$\frac{d\|{\bf
x}\|}{dt} =    \sum_{j}\Big(\sum_{i} h_{ij}\Big)x_{j} \ge h_*(t)
\sum_{j} x_j = h_*(t) \|{\bf x}\|,$$ and $\|{\bf x}(t)\| \ge
e^{\int_0^t h_*(\tau)\,d\tau}\|{\bf x}(0)\|$.

\begin{remark} If we can find $D$ such that $h^*(t)= h_*(t)$, then
Theorem 1 gives {\it sharp} bound on the rate of convergence.
\end{remark}

\begin{remark} Let $X(t)$ be a homogeneous Markov chain. Then the
corresponding spectral gap, or decay parameter,
$$ \beta := \sup \{a>0: \|{\bf p}(t)-{\bf \pi}\| = \mathcal{O}(e^{-at}) ~~\mbox{as}~~
t\to\infty ~~\mbox{for~all}~~{\bf p}(0)\},$$ satisfies the
inequality $h_* \le \beta \le  h^*$. This bound makes it possible to find the
asymptotic behavior of the spectral gap as the number of states
tends to infinity, see, for instance $($Granovsky and Zeifman, 1997,
2000, 2005; Van Doorn et al, 2010$)$.
\end{remark}

\section{Birth-death-catastrophe  process}

Let $X(t)$ be a birth-death-catastrophe  process (BDPC) on the finite
state space $\left\{0,\dots,S\right\}$ and let $ \lambda_n\left(
t\right)$, $\mu_{n+1}\left( t\right)$, and $\xi_{n+1}(t), \
n=0,\dots,S-1$ be the corresponding birth,  death and catastrophe
intensities.  Then
\begin{equation}
a_{ij} (t) = \left\{
\begin{array}{ll}
\lambda_{i-1} \left( t\right), & \mbox {\it if } j=i-1, \\
\mu_{i+1}\left( t\right), & \mbox {\it if } j=i+1>1,
\\
-\left( \lambda_i\left( t\right) +\mu_i \left( t\right) + \xi_i(t)\right), & \mbox {\it if } j=i, \\
\xi_j(t), & \mbox {\it if } i=0,\,j>1, \\
\mu_1(t)+\xi_1(t), & \mbox {\it if } i=0,\, j=1, \\
0, &  \mbox{\it overwise}
\end{array}
\right.  \label{401}
\end{equation}
are the transposed intensities of the process. A detailed study of
this class of processes can be found in (Zeifman et al, 2013). Then
the assumptions (i) and (ii) of  Section 3 are fulfilled for the upper
triangular matrix
\begin{equation}
  D=\left(
  \begin{array}{ccccc}
  d_1 & d_1 & d_1 & \cdots & d_1 \\
  0   & d_2 & d_2 & \cdots & d_2 \\
  0   & 0   & d_3 & \cdots & d_3 \\
  \vdots & \vdots & \vdots & \ddots \\
  0 & 0 & 0 & 0 & d_{S}
  \end{array}
  \right).
\label{402}
\end{equation}

Let $\alpha _k(t)=-\sum_{i} h_{ik}(t)$. Then
\begin{equation}
\alpha _k(t)=\lambda_k(t)+\mu_{k+1}(t)+ \xi_{k+1}(t) -\frac{d_{k+1}}{d_k} \lambda_{k+1}(t)-\frac{d_{k-1}}{d_k} \mu_k(t), \quad k=0, \dots, S-1, \label{403}
\end{equation}
and
\begin{equation}
h^*(t) = -\beta^*(t) = - \min\limits_{0 \le k \le S-1} \alpha _k(t), \quad h_*(t) = - \beta_*(t) = -\max\limits_{0 \le k \le S-1} \alpha _k(t).\label{404}
\end{equation}
Hence, the following statement holds.

\begin{theorem} Let $X(t)$ be a finite BDPC, and let there exist  positive numbers $\{d_i\}$ such that
 \begin{equation}
 \int_0^\infty \beta^*(t)\, dt = +\infty.
 \label{405}
\end{equation}
Then $X(t)$ is weakly ergodic and the following bounds hold:
\begin{equation}
\|{\bf z^{*}}(t) - {\bf z^{**}}(t)\|_{1D} \le e^{-\int_0^t
\beta^*(\tau)\,d\tau}\|{\bf z^{*}}(0) - {\bf z^{**}}(0)\|_{1D},
\label{406}
\end{equation}
\noindent for any corresponding initial conditions ${\bf p^{*}}(0)$,
${\bf p^{**}}(0)$, and
\begin{equation}
\|{\bf z^{*}}(t) - {\bf z^{**}}(t)\|_{1D} \ge e^{-\int_0^t \beta_*(\tau)\,d\tau}\|{\bf z^{*}}(0) - {\bf z^{**}}(0)\|_{1D},
\label{407}
\end{equation}
\noindent if the initial conditions are such that $D\left({\bf z^{*}}(0)
- {\bf z^{**}}(0)\right) \ge {\bf 0}$.
\end{theorem}

\begin{example} Let $ \lambda_n\left(t\right) = \lambda(t)$,
$\mu_{n+1}\left( t\right) =(n+1)\mu(t)$, for $n=0,\dots,S-1$,
$\xi_{S}(t) \equiv 0$, and $\xi_{n+1}(t) =\lambda(t)$, if
$n=0,\dots,S-2$. Put all $d_i =1$. Then in $($\ref{403}$)$ all $\alpha _k(t)=\lambda(t)
+ \mu(t)$, and Theorem 2 gives the inequality
\begin{equation}
\|{\bf z^{*}}(t) - {\bf z^{**}}(t)\|_{1D} \le e^{-\int_0^t
\left(\lambda(\tau)+\mu(\tau)\right)\,d\tau}\|{\bf z^{*}}(0) - {\bf z^{**}}(0)\|_{1D},
\label{408}
\end{equation}
\noindent for any initial conditions. Moreover, if the initial conditions are such
that $D\left({\bf z^{*}}(0) - {\bf z^{**}}(0)\right) \ge {\bf 0}$, then
\begin{equation}
 \|{\bf z^{*}}(t) - {\bf z^{**}}(t)\|_{1D} = e^{-\int_0^t
\left(\lambda(\tau)+\mu(\tau)\right)\,d\tau}\|{\bf z^{*}}(0) - {\bf z^{**}}(0)\|_{1D}.
\label{409}
\end{equation}
\end{example}

\begin{remark} A finite birth-death process with constant rates of
birth $\lambda_k(t) = a$ and $\mu_{k+1}(t)=b$ was considered in
$($Granovsky and Zeifman. 1997$)$, where the corresponding $D$ and
sharp $\beta_*=\beta^*= a+b-2\sqrt{ab}\cos\frac{\pi}{S+1} \to
(\sqrt{a}-\sqrt{b})^2$ as $S \to \infty$ were obtained.
\end{remark}

\section{SZK model}

Let $X=X(t)$, $t\geq 0$, be a Markov chain with intensities
$q_{i,i+k}\left( t\right) =\lambda_k(t)$, $q_{i,i-k}\left( t\right)
= \mu_k(t)$ for any $k > 0$ (SZK model, see (Satin et al, 2013)). In other words, in terms of
queueing theory, we suppose that the arrival rates $\lambda_{k}(t)$
and the service rates $\mu_{k}(t)$ do not depend on the the length
of the queue. In addition, we assume that $\lambda_{k+1}(t) \le
\lambda_{k}(t)$ and  $\mu_{k+1}(t) \le \mu_{k}(t)$ for any $k$ and
almost all $t \ge 0$. Then
%{\scriptsize
\begin{equation}
A(t)=\left(
\begin{array}{cccccccc}
a_{00}(t) & \mu_1(t)  & \mu_2(t)   & \mu_3(t)  &  \cdots & \mu_S(t) \\
\la_1(t)   & a_{11}(t)  & \mu_1(t)  & \mu_2(t)     & \cdots & \mu_{S-1}(t) \\
\la_2(t)  & \la_1(t)    & a_{22}(t)& \mu_1(t)      &  \cdots & \mu_{S-2}(t) \\
\cdots \\
\la_S(t) & \la_{S-1}(t) & \la_{S-2}(t) & \cdots & \cdots    &  a_{SS}(t)\\
\end{array}
\right),
\label{501}
\end{equation}
%}
\noindent where  $a_{ii}(t)=-\sum_{k=1}^{i}\mu_k(t) -
\sum_{k=1}^{S-i} \la_{k}(t)$ is the transposed intensity matrix of
the process. This class of processes was introduced and studied in
our previous papers (Satin et al, 2013, Zeifman et al, 2014b). Then
the assumptions (i) and (ii) of  Section 3 are fulfilled for upper
triangular matrix (\ref{402}). In this case,
$$
D BD^{-1} =\hspace{10cm}
$$
\begin{equation}
\left(\!\!
\begin{array}{ccccccc}
a_{11}-\la_S  &(\mu_1-\mu_2)\frac{d_1}{d_2}  & (\mu_2-\mu_3)\frac{d_1}{d_3}  & \cdots & (\mu_{S-1}-\mu_S)\frac{d_1}{d_S} \vspace{2mm}\\
(\la_1-\la_S)\frac{d_2}{d_1} &  a_{22}-\la_{S-1}  & (\mu_1-\mu_3)\frac{d_2}{d_3}  & \cdots &  (\mu_{S-2}-\mu_S)\frac{d_2}{d_S}\vspace{2mm} \\
(\la_2-\la_S) \frac{d_3}{d_1} &  (\la_1-\la_{S-1}) \frac{d_3}{d_2}  &a_{33}-\la_{S-2} & \cdots & (\mu_{S-3}-\mu_S) \frac{d_3}{d_S}\vspace{2mm}\\
\ldots&\ldots&\ldots&\ldots&\ldots\vspace{2mm}\\
(\la_{S-1} -\la_S) \frac{d_S}{d_1} & (\la_{S-2} -\la_{S-1}) \frac{d_S}{d_2}  & (\la_{S-3} -\la_{S-2})\, \frac{d_S}{d_3}  & \cdots & a_{SS}-\la_1
\end{array}
\!\right),
\end{equation}
\noindent and $-\alpha _k(t) = \sum_{i} h_{ik}(t)$ equals to   the sum of elements of the $k-$th column of $D B(t)D^{-1}$.
Hence, the following statement holds.

\begin{theorem} Let $X(t)$ be a finite SZK model, and let there exist  positive numbers $\{d_i\}$ such that
 \begin{equation}
 \int_0^\infty \beta^*(t)\, dt = +\infty.
 \label{505}
\end{equation}
Then $X(t)$ is weakly ergodic and the following bounds hold:
\begin{equation}
\|{\bf z^{*}}(t) - {\bf z^{**}}(t)\|_{1D} \le e^{-\int_0^t
\beta^*(\tau)\,d\tau}\|{\bf z^{*}}(0) - {\bf z^{**}}(0)\|_{1D},
\label{506}
\end{equation}
\noindent for any corresponding initial conditions ${\bf p^{*}}(0)$,
${\bf p^{**}}(0)$, and
\begin{equation}
\|{\bf z^{*}}(t) - {\bf z^{**}}(t)\|_{1D} \ge e^{-\int_0^t \beta_*(\tau)\,d\tau}\|{\bf z^{*}}(0) - {\bf z^{**}}(0)\|_{1D},
\label{507}
\end{equation}
\noindent if the initial conditions are such that $D\left({\bf z^{*}}(0)
- {\bf z^{**}}(0)\right) \ge {\bf 0}$.
\end{theorem}

\begin{example} Let  $\lambda_{k}(t) = \frac{\lambda(t)}{k}$,
$\mu_{k} =\frac{\mu}{k}$. Let $X(t)$ be a queue-length process for
the corresponding queueing model with  batch arrivals and group
services. Let $d_k=1$, $k\ge1$. Then $\alpha _k(t)=\lambda(t) +
\mu(t)$, $k\ge1$, and Theorem 3 gives the inequality
\begin{equation}
\|{\bf z^{*}}(t) - {\bf z^{**}}(t)\|_{1D} \le e^{-\int_0^t
\left(\lambda(\tau)+\mu(\tau)\right)\,d\tau}\|{\bf z^{*}}(0) - {\bf
z^{**}}(0)\|_{1D}, \label{508}
\end{equation}
\noindent for any initial conditions. Moreover, if the initial
conditions are such that $D\left({\bf z^{*}}(0) - {\bf z^{**}}(0)\right)
\ge {\bf 0}$, then
\begin{equation}
 \|{\bf z^{*}}(t) - {\bf z^{**}}(t)\|_{1D} = e^{-\int_0^t
\left(\lambda(\tau)+\mu(\tau)\right)\,d\tau}\|{\bf z^{*}}(0) - {\bf
z^{**}}(0)\|_{1D}. \label{509}
\end{equation}
\end{example}

\section{Markov chain with absorbtion in zero}

Let now $X(t)$, $t\geq 0$ be a Markov chain with absorption in zero,
i. e.,  let $q_{00}(t) \equiv 0$. Then
\begin{equation}
{%\scriptsize
B = \left(
\begin{array}{ccccccccc}
a_{11}   & a_{12}   &  \cdots & a_{1S}  \\
a_{21}  & a_{22}    &    \cdots & a_{2S}  \\
a_{31}     & a_{32}   &    \cdots & a_{3S}  \\
\cdots \\
a_{S1}  & a_{S2} & \cdots     &  a_{SS}
\end{array}
\right)}\label{601}
\end{equation}
\noindent is essentially nonnegative itself. Hence,  the assumptions
(i) and (ii) of  Section 3 are fulfilled for the diagonal matrix
\begin{equation}
  D=\left(
  \begin{array}{ccccc}
  d_1 & 0 & 0 & \cdots & 0 \\
  0   & d_2 & 0 & \cdots & 0 \\
  0   & 0   & d_3 & \cdots & 0 \\
  \vdots & \vdots & \vdots & \ddots \\
  0 & 0 & 0 & 0 & d_{S}
  \end{array}
  \right).
\label{602}
\end{equation}

If $\alpha _k(t)=-\sum_{i} h_{ik}(t)$, then
\begin{equation}
\alpha _k(t)=a_{kk}(t) - \sum_{i\neq k} \frac{d_{i}}{d_k} a_{ik}(t),
\quad k=1, \dots, S, \label{603}
\end{equation}
and
\begin{equation}
h^*(t) = -\beta^*(t) = - \min\limits_{1 \le k \le S} \alpha _k(t),
\quad h_*(t) = - \beta_*(t) = -\max\limits_{1 \le k \le S} \alpha
_k(t).\label{604}
\end{equation}
Hence, the following statement holds.

\begin{theorem} Let $X(t)$ be a finite a Markov chain with absorption in zero, and let there exist  positive numbers $\{d_i\}$ such that
 \begin{equation}
 \int_0^\infty \beta^*(t)\, dt = +\infty.
 \label{605}
\end{equation}
Then $X(t)$ is weakly ergodic and the following bounds hold:
\begin{equation}
\|{\bf z^{*}}(t) - {\bf z^{**}}(t)\|_{1D} \le e^{-\int_0^t
\beta^*(\tau)\,d\tau}\|{\bf z^{*}}(0) - {\bf z^{**}}(0)\|_{1D},
\label{606}
\end{equation}
\noindent for any corresponding initial conditions ${\bf p^{*}}(0)$,
${\bf p^{**}}(0)$, and
\begin{equation}
\|{\bf z^{*}}(t) - {\bf z^{**}}(t)\|_{1D} \ge e^{-\int_0^t
\beta_*(\tau)\,d\tau}\|{\bf z^{*}}(0) - {\bf z^{**}}(0)\|_{1D},
\label{607}
\end{equation}
\noindent if the initial conditions are such that ${\bf z^{*}}(0) - {\bf
z^{**}}(0) \ge {\bf 0}$.
\end{theorem}

\begin{example} Let now $X(t)$ be a BDP with birth rates
$\lambda_0(t) \equiv 0$, $\lambda_k(t)=2\phi(t)$, $1 \le k <S$,
and death rates $\mu_1(t)=3\phi(t)$, $\mu_k(t)=6\phi(t)$, $1 < k<S$,
$\mu_S(t)=2\phi(t)$. Let $d_k=2^{k-1}$, $1\le k\le S$. Then $\alpha _k(t)=\phi(t)$, $1\le k\le S$,
and Theorem 4 gives the inequality
\begin{equation}
\|{\bf z^{*}}(t) - {\bf z^{**}}(t)\|_{1D} \le e^{-\int_0^t
\phi(\tau)\,d\tau}\|{\bf z^{*}}(0) - {\bf z^{**}}(0)\|_{1D},
\label{608}
\end{equation}
\noindent for any initial conditions. Moreover, if the initial
conditions are such that ${\bf z^{*}}(0) - {\bf z^{**}}(0) \ge {\bf 0}$,
that
\begin{equation}
 \|{\bf z^{*}}(t) - {\bf z^{**}}(t)\|_{1D} = e^{-\int_0^t
\phi(\tau)\,d\tau}\|{\bf z^{*}}(0) - {\bf z^{**}}(0)\|_{1D}.
\label{609}
\end{equation}
\end{example}

\bigskip

{\bf Acknowledgement.} This work was supported by the Russian Scientific
Foundation  (Grant No. 14-11-00397).

\newpage

\begin{center} {\bf References} \end{center}

\noindent Van Doorn~E.\,A., Zeifman~A.\,I.,  Panfilova~T.\,L. 2010.
Bounds and asymptotics for the rate of convergence of birth-death
processes.  Th. Prob. Appl., 54, 97--113.

\medskip

\noindent Granovsky, B. L.,  Zeifman, A. I. 1997. The decay function
of nonhomogeneous birth-death processes, with application to
mean-field models. Stochastic processes and their applications, 72,
105--120.

\medskip

\noindent Granovsky~B.\,L.,  Zeifman~A.\,I. 2000.  The
$\footnotesize N$-limit of spectral gap of a class of birth-death
Markov chains.  Appl.\ Stoch.\ Models Bus.\ Ind., 16, No 4,
235--248.

\medskip

\noindent Granovsky, B.~L., Zeifman, ~ A.~I. 2004. Nonstationary
Queues: Estimation of the Rate of Convergence. Queueing Syst.  46,
363--388.

\medskip

\noindent Granovsky, B. L.,  Zeifman, A. I. 2005. On the lower bound
of the spectrum of some mean-field models. Theory of Probability \&
Its Applications, 49, 148--155.

\medskip

\noindent Kartashov, N.~V. 1985. Criteria for uniform ergodicity and
strong stability of Markov chains with a common phase space, Theory
Probab. Appl., 30, 71 -- 89.

\medskip

\noindent  Kartashov, N.~V. 1996.  Strong stable Markov chains,
Utrecht: VSP. Kiev: TBiMC.

\medskip

\noindent  Liu, Y. 2012.  Perturbation bounds for the stationary
distributions of Markov chains. SIAM. J. Matrix Anal. \& Appl., 33,
1057--1074.

\medskip

\noindent Mitrophanov, A.~Yu. 2003. Stability and exponential
convergence of continuous-time Markov chains. J. Appl. Prob.  40,
970--979.

\medskip

\noindent Mitrophanov, A.~Yu. 2004. The spectral gap and
perturbation bounds for reversible continuous-time Markov chains. J.
Appl. Probab. 41, 1219--1222.

\medskip

\noindent Satin, Ya.~A., Zeifman, A.~I., Korotysheva, ~A.~V. 2013.
On the rate of convergence and truncations for a class of Markovian
queueing systems. Theory Probab. Appl., 57, 529--539.

\medskip

\noindent  Zeifman, A.~I.  1985.  Stability for contionuous-time
nonhomogeneous Markov chains,  Lecture Notes in Mathematics, 1233.
Berlin etc.: Springer-Verlag, 401--414.

\medskip

\noindent  Zeifman, A.~I.  1995. Upper and lower bounds on the rate
of convergence for nonhomogeneous birth and death processes. Stoch.
Proc. Appl., 59, 157--173.

\medskip

\noindent Zeifman, A.~I.,  Leorato,  S., Orsingher, E., Satin, Ya.,
Shilova, G. 2006. Some universal limits for nonhomogeneous birth and
death processes. Queueing Syst.  52, 139--151.

\medskip

\noindent Zeifman, A.~I., Satin, Ya.,  Panfilova, T.. 2013 Limiting
characteristics for finite birth-death-catastrophe processes.
Mathematical Biosciences, 245, 96--102.

\medskip

\noindent  Zeifman, A.~I.,  Korolev, V.~Y. 2014a. On perturbation
bounds for continuous-time Markov chains. Statistics \& Probability
Letters, 88, 66--72.

\medskip

\noindent Zeifman, A.~I., Korolev, V., Korotysheva, A., Satin, Y.,
Bening, V. 2014b. Perturbation bounds and truncations for a class of
Markovian queues. Queueing Systems, 76, 205--221.

\end{document}